\documentclass[12pt]{article}
\usepackage{bbm}
\usepackage{amssymb,amsmath}

\newtheorem{theorem}{Theorem}[section]
\newtheorem{lemma}[theorem]{Lemma}

\newtheorem{definition}[theorem]{Definition}
\newtheorem{corollary}[theorem]{Corollary}
\newtheorem{example}[theorem]{Example}

\begin{document}

\title{Admissible operators and ${\mathcal H}_{\infty}$ calculus
}

\author{Hans Zwart\thanks{University of Twente, Faculty of Electrical Engineering, Mathematics and Computer Science, Department of Applied Mathematics, P.O. Box 217, 7500 AE Enschede, The Netherlands, \tt h.j.zwart@math.utwente.nl}
}

\maketitle

\begin{abstract}
  Given a Hilbert space and the generator $A$ of a strongly
  continuous, exponentially stable, semigroup on this Hilbert space.
  For any $g(-s) \in {\mathcal H}_{\infty}$ we show that there exists
  an infinite-time admissible output operator $g(A)$. If $g$ is
  rational, then this operator is bounded, and equals the ``normal''
  definition of $g(A)$. In particular, when $g(s)=1/(s + \alpha)$, $
  \alpha \in {\mathbb C}_0^+$, then this admissible output operator
  equals $(\alpha I - A)^{-1}$.

  Although in general $g(A)$ may be unbounded, we always have that
  $g(A)$ multiplied by the semigroup is a bounded operator for every
  (strictly) positive time instant. Furthermore, when there exists an
  admissible output operator $C$ such that $(C,A)$ is exactly
  observable, then $g(A)$ is bounded for all $g$'s with $g(-s) \in
  {\mathcal H}_{\infty}$, i.e., there exists a bounded ${\mathcal H}_{\infty}$-calculus. Moreover, we rediscover some well-known classes of generators also having a bounded ${\mathcal H}_{\infty}$-calculus.
\end{abstract}


\section{Introduction}

Functional calculus is a sub-field of mathematics with a long history.
It started in the thirties of the last century with the work by von
Neumann for self-adjoint operators \cite{Neum96}, and was further
extended by many researchers, see e.g.\ \cite{HiPh57} and
\cite{DuSc71}. For an overview, see the book by Markus Haase,
\cite{Haas06}. The basic idea behind functional calculus for the
operator $A$ is to construct a mapping from an algebra of (scalar)
functions to the class of (bounded) operators, such that
\begin{itemize}
\item The function identically equals to one is mapped to the identity
  operator;
\item If $f(s) = (s-a)^{-1}$, then $f(A) = (sI-A)^{-1}$;
\item Furthermore, the operator associated to $f_1 \cdot f_2$ equals
  $f(A) f_2(A)$.
\end{itemize} 

Before we explain the contribution of this paper, we introduce some
notation. By $X$ we denote separable Hilbert space with inner product
$\langle \cdot,\cdot \rangle$ and norm $\|\cdot\|$, and by $A$ we
denote an unbounded operator from its domain $D(A) \subset X$ to $X$.
We assume that $A$ generates an exponentially stable semigroup on $X$,
which we denote by $\left( T(t) \right)_{t \geq 0}$.

By ${\mathcal H}_{\infty}^-$ we denote the space of all bounded,
analytic functions defined on the half-plane ${\mathbb C}^- := \{ s
\in {\mathbb C} \mid \text{Re}(s) < 0\}$. It is clear that this
function class is an algebra under pointwise multiplication and
addition. Hence this could serve as a class for which one could build
a functional calculus. However, it is known that there exists a
generator of exponential stable semigroup, which does not have a
functional calculus with respect to ${\mathcal H}_{\infty}^-$. For
proof of this and many more we refer to \cite{AlDM96}, \cite{Haas06},
and the references therein. Although a bounded functional calculus is
not possible, an unbounded functional calculus is always possible.
\begin{theorem}
  \label{T:0}%
  Under the assumptions stated above, we have that for all $g \in
  {\mathcal H}_{\infty}^-$ there exists an operator $g(A)$ which is
  bounded from the domain of $A$ to $X$, and which is admissible,
  i.e.,
  \[
    \int_0^{\infty} \|g(A) T(t) x_0 \|^2 dt \leq \gamma_A \|g\|_{\infty}^2 \|x_0\|^2, \qquad x_0 \in X.
  \]
  The mapping $g \mapsto g(A)$ satisfies the conditions of a
  functional calculus. Furthermore, for all $t >0$, we have that $g(A)
  T(t)$ can be extended to a bounded operator, and
  \[
    \|g(A) T(t) \| \leq \frac{\gamma }{\sqrt{t}}.
  \]
\end{theorem}

Apart from proving this theorem, we shall also rediscover some classes
of generators for which $g(A)$ is bounded for all $g \in {\mathcal
  H}_{\infty}^-$, i.e., for which there is a bounded functional
calculus.

For the proof of the above result, we need beside the Hardy space
${\mathcal H}_{\infty}^-$ also the Hardy spaces ${\mathcal H}_2(X)$
and ${\mathcal H}_2^{\perp}(X)$.

${\mathcal H}_2(X)$ and ${\mathcal H}_2^{\perp}(X)$ denote the Laplace
transform, ${\mathcal L}$, of functions in $L^2((0,\infty),X)$ and
$L^2((-\infty,0),X)$, respectively. It is known that this
transformation is an isometry.  Every function in ${\mathcal
  H}_{\infty}^-, {\mathcal H}_2(X)$ and ${\mathcal H}_2^{\perp}(X)$
has a unique extension to the imaginary axis on which this functions are
bounded, and square integrable, respectively. Furthermore, the norm of
$g \in {\mathcal H}_{\infty}^-$ equals the (essential) supremum over
the imaginary axis of the boundary function. Let $f(t)$ be a function
in $L^2((0,\infty),X)$ with Laplace transform $F(s)$, and let
$f_{\mathrm{ext}}(t)$ be the function in $L^2((-\infty,\infty),X)$
defined by
\[
  f_{\mathrm{ext}}(t) = \begin{cases}
  f(t) & t \geq 0\\ 0 & t < 0
  \end{cases}
\]
Then the Fourier transform $\hat{f}_{\mathrm{ext}}$ of $f_{\mathrm{ext}}(t)$ satisfies $\hat{f}_{\mathrm{ext}}(\omega) = F(i\omega)$, for almost all $\omega \in {\mathbb R}$. Here $F(i \cdot)$ denote the boundary function of the Laplace transform $F(s)$.

We define the following Toeplitz operator on $L^2((0,\infty);X)$
\begin{definition}
  \label{D1.1}
  Let $g$ be an element of ${\mathcal H}_{\infty}^-$. Associated to
  this function we define the mapping $M_g$ as
  \begin{equation}
  \label{eq:1}
  M_g f = {\mathcal L}^{-1} \left(\Pi\left( g F\right)\right), \qquad f \in L^2((0,\infty), X),
  \end{equation}
  where $F$ denotes the Laplace transform of $f$. $\Pi$ denotes the
projection onto ${\mathcal H}_2(X)$.
\end{definition}

It is clear that this is a linear bounded map from $L^2((0,\infty);X)$
into itself, and
\begin{equation}
  \label{eq:1.1}
  \|M_g \| \leq \|g\|_{\infty}.
\end{equation}
Furthermore, it follows easily from (\ref{eq:1}) that if
$K$ is a bounded mapping on $X$, then its commutes with $M_g$, i.e.,
\begin{equation}
  \label{eq:7}
  K M_g = M_g K.
\end{equation}
 
It is easy to see that ${\mathcal H}_{\infty}^-$ is an algebra under the multiplication and addition. In particular $g_1 g_2 \in {\mathcal H}_{\infty}^-$ whenever $g_1$, $g_2 \in {\mathcal H}_{\infty}^-$. Furthermore, we have the following result.
\begin{lemma}
  \label{L1.2}
  Let $g_1$ and $g_2$ be elements of ${\mathcal H}_{\infty}^-$. Then
  \begin{equation}
    \label{eq:1.2}
    M_{g_1g_2} = M_{g_1} M_{g_2}.
  \end{equation}
  In particular, if $g$ is invertible in ${\mathcal H}_{\infty}^-$, then $M_g$ is (boundedly) invertible and $\left(M_g\right)^{-1} = M_{g^{-1}}$.
\end{lemma}
{\bf Proof}\/
We use the fact that any $g \in {\mathcal H}_{\infty}^-$ maps ${\mathcal H}_2^{\perp}$ into ${\mathcal H}_2^{\perp}$.
\begin{eqnarray*}
  M_{g_1} M_{g_2} f &=& {\mathcal L}^{-1} \left( \Pi g_1 \left( \Pi\left( g_2 F\right)\right)\right) \\
  &=&
  {\mathcal L}^{-1} \left( \Pi \left(g_1 g_2 F\right)\right) +  {\mathcal L}^{-1} \left( \Pi \left( g_1 (I-\Pi)\left( g_2 F\right)\right)\right)\\
  &=& {\mathcal L}^{-1} \left( \Pi \left(g_1 g_2 F\right)\right) +0,
\end{eqnarray*}
where we have used the above mentioned fact that $ g_1 (I-\Pi)$ maps
into ${\mathcal H}_2^{\perp}$, and so $\Pi g_1 (I-\Pi)=0$. Since by
definition ${\mathcal L}^{-1} \left( \Pi \left(g_1 g_2
    F\right)\right)$ equals $M_{g_1g_2} f$, we have proved the first
assertion.

The last assertion follows directly, since $M_1 = I$.
\hfill$\Box$
\medskip

By $\sigma_{\tau}$ we denote the shift with $\tau\geq 0$, i.e.,
\begin{equation}
  \label{eq:2}
  \left(\sigma_\tau (f)\right)(t) = f(t+\tau),\qquad t \geq 0.
\end{equation}
This is also a linear bounded map from $L^2((0,\infty);X)$ into
itself. This mapping commutes with $M_g$ as is shown next.
\begin{lemma}
\label{L1.3}
  For all $\tau >0$ and all $g$ in ${\mathcal H}_{\infty}^-$, we have that
  \begin{equation}
  \label{eq:3}
  \sigma_{\tau} \left( M_g f \right) = M_g \left(\sigma_{\tau} f\right), \qquad f \in L^2((0,\infty), X).
\end{equation}
\end{lemma}
{\bf Proof}\/ We use the following well-known equality. If $h$ is Fourier transformable, then the Fourier transform of $h(\cdot + \tau)$ equals $e^{i \omega \tau} \hat{h}(\omega)$, where $\hat{h}$ denotes the Fourier transform of $h$.

Let $h \in L^2((0,\infty);X)$, then 
\begin{equation}
\label{eq:1.3}
  {\mathcal L}(\sigma_\tau h) = \widehat{\left(\sigma_{\tau} h\right)_{\mathrm{ext}}} = \widehat{\sigma_{\tau} h_{\mathrm{ext}}} - \hat{q} = e^{i\omega \tau} \widehat{h_{\mathrm{ext}}} - \hat{q} = e^{i\omega \tau} {\mathcal L}(h) - \hat{q},
\end{equation}
with $q \in L^2((-\infty,0);X)$.
In particular, we find for every $h \in L^2(0,\infty);X)$ that
\begin{equation}
  \label{eq:1.4}
  {\mathcal L}(\sigma_\tau h) = \Pi \left({\mathcal L}(\sigma_\tau h)\right) = \Pi \left(e^{i\omega \tau} {\mathcal L}(h) \right) - 0 ={\mathcal L}\left( M_{e^{i\cdot \tau}} h\right),
\end{equation}
where we have used that $e^{i\omega \tau}$ is the boundary function corresponding to $e^{i s \tau} \in {\mathcal H}_{\infty}^-$.

Using (\ref{eq:1.3}) we see that
\begin{equation}
\label{eq:1.5}
  M_g \left(\sigma_{\tau} f \right) = {\mathcal L}^{-1} \left(\Pi\left( g e^{i\cdot \tau} {\mathcal L}(f)\right) \right) - {\mathcal L}^{-1} \left(\Pi \left( g \hat{q} \right)\right) =
  {\mathcal L}^{-1} \left(\Pi\left( g e^{i\cdot \tau} {\mathcal L}(f)\right) \right),
\end{equation}
since $\hat{q} \in {\mathcal H}_2^{\perp}(X)$, and since $g \in {\mathcal H}_{\infty}^-$.
Using Lemma \ref{L1.2}, we find that
\begin{equation}
  \label{eq:1.6}
  M_g \left(\sigma_{\tau} f \right) = {\mathcal L}^{-1} \left(\Pi\left( g e^{i\cdot \tau} {\mathcal L}(f)\right) \right) = M_{e^{i\cdot \tau} g}f = M_{e^{i\cdot \tau}} M_g f.
\end{equation}
Now using (\ref{eq:1.4}), we see that 
\begin{equation}
  \label{eq:1.7}
  M_g  \left(\sigma_{\tau} f \right) = \sigma_{\tau} \left( M_g f \right).
\end{equation}
\hfill$\Box$

\section{Output maps and admissible output operators}

In this section we study admissible operators which commute with the semigroup. We begin by defining well-posed output maps.
\begin{definition}
  \label{D:2.0}  Let $\left( T(t)\right)_{t\geq 0}$ be a strongly continuous semigroup on the Hilbert space $X$, and let $Y$ be another Hilbert space. We say that the mapping ${\mathcal O}$ is a well-posed (infinite-time) output map if
  \begin{itemize}
  \item ${\mathcal O}$ is a bounded linear mapping from $X$ into $L^2((0,\infty);Y)$, and
  \item For all $\tau \geq 0$ and all $x_0 \in X$, we have that $\sigma_{\tau} {\mathcal O} x_0 = {\mathcal O}\left(T(\tau) x_0\right)$.
  \end{itemize}
\end{definition}

Closely related to well-posed output mappings are admissible operators, which are defined next.
\begin{definition}
  \label{Dn:2.1}%
  Let $\left( T(t)\right)_{t\geq 0}$ be a strongly continuous
  semigroup on the Hilbert space $X$. Let $D(A)$ be the domain of its
  generator $A$. A linear mapping $C$ from $D(A)$ to $Y$, another
  Hilbert space, is said to be an (infinite-time) admissible output
  operator for $\left( T(t)\right)_{t\geq 0}$ if $CT(\cdot) x_0 \in
  L^2((0,\infty),Y)$ for all $x_0 \in D(A)$ and there exists an $m$
  independent of $x_0$ such that
  \begin{equation}
  \label{eq:2.1}
    \int_0^{\infty} \| C T(t) x_0 \|_Y^2 dt \leq m \|x_0\|_X^2.
  \end{equation}
\end{definition}

If $C$ is (infinite-time) admissible, then for all $x_0 \in X$ we can
uniquely define an $L^2((0,\infty),X)$-function. We denote this
function by $C T(\cdot) x_0$. Hence ${\mathcal O} : X \rightarrow
L^2((0,\infty);Y)$ defined by ${\mathcal O}x_0 = C T(\cdot) x_0$ is a
well-posed output map. From \cite{Weis89} we know that the converse
holds as well.
\begin{lemma}
  \label{L:2.1a} If ${\mathcal O}$ is a well-posed output mapping, then there exists a (unique) linear bounded mapping from $D(A)$ to $Y$, $C$, such that ${\mathcal O}x_0 = C T(\cdot) x_0$ for all $x_0$.
\end{lemma}

In the sequel of this section we concentrate on admissible output
operators which commute with the semigroup, i.e., $C$ a linear operator from
$D(A)$ to $X$ and
\begin{equation}
  \label{eq:2.2}
  C T(t)x_0 = T(t) C x_0 \qquad \mbox{ for all } t \geq 0 \mbox{ and } x_0 \in D(A).
\end{equation}
For these operators we have the following results.
\begin{lemma}
  \label{L:2.1b} Let $C$ be the admissible output operator associated with the well-posed output map ${\mathcal O}$. Then (\ref{eq:2.2}) holds if and only if for all $t \geq 0$ there holds ${\mathcal O} T(t) = T(t){\mathcal O}$.
\end{lemma}
\begin{theorem}
  \label{T2.2} Let $C$ be a bounded linear operator from $D(A)$ to $X$, which is admissible for the exponentially stable semigroup $\left( T(t) \right)_{t \geq 0}$ and which commutes with this semigroup. Then the following holds
    \begin{enumerate}
    \item For all $x_0 \in D(A)$, we have that $C A^{-1} x_0 = A^{-1} C x_0$.
    \item For all $t >0$, the operator $CT(t): D(A) \rightarrow X$ can
      be extended to a bounded operator on $X$. Furthermore, $\|C
      T(t)\| \leq \gamma t^{-1/2}$ for some $\gamma$ independent of
      $t$.
    \end{enumerate}
\end{theorem}
{\bf Proof}\/ The first assertion follows easily from (\ref{eq:2.2}) by using Laplace transforms. We concentrate on the second assertion.

Let $x_0\in D(A)$ and $x_1 \in X$, then for $t >0$ we have that
\begin{eqnarray*}
  t \langle x_1, C T(t) x_0 \rangle &=& \int_0^t \langle x_1, C T(t) x_0 \rangle d \tau \\
  &=&
  \int_0^t \langle x_1, C T(\tau) T(t-\tau) x_0 \rangle d \tau \\
  &=&
   \int_0^t \langle x_1,  T(\tau) C T(t-\tau) x_0 \rangle d \tau \\
  &=&
  \int_0^t \langle T(\tau)^* x_1, C T(t-\tau) x_0 \rangle d\tau \\
  &\leq& \sqrt{\int_0^t \|T(\tau)^* x_1 \|^2 d\tau} \sqrt{\int_0^t \|CT(t-\tau) x_0 \|^2 d\tau}.
\end{eqnarray*}
Using the fact that the semigroup, and hence its adjoint, are uniformly bounded, and the fact that $C$ is (infinite-time) admissible, we find that
\[
  t \langle x_1, C T(t) x_0 \rangle \leq \sqrt{t} M \|x_1\| m \|x_0\|.
\]
Since this holds for all $x_1 \in X$, we conclude that
\[
  t \|CT(t) x_0 \| \leq \sqrt{t}  mM \|x_0\|.
\]
This inequality holds for all $x_0 \in D(A)$. The domain of a generator is dense, and hence we have proved the second assertion.
\hfill$\Box$
\medskip

\mbox{}From Theorem \ref{T2.2} it is clear that if the semigroup is surjective, then any admissible $C$ which commutes with the semigroup is bounded. However, this does not hold for a general semigroup as is shown in the following example. Furthermore, this example also shows that the estimate in the previous theorem cannot be improved.
\begin{example}
\label{E2.6}
 Let $\{\phi_n, n \in {\mathbb N}\}$ be an orthonormal basis of $X$, and define for $t \geq 0$ the operator
 \begin{equation}
   \label{eq:9}
    T(t) \sum_{n=1} ^N \alpha_n \phi_n =  \sum_{n=1} ^N  e^{-n^2t}  \alpha_n \phi_n .
\end{equation}
 It is not hard to show that this defines an exponentially stable $C_0$-semigroup on $X$. The infinitesimal generator $A$ is given by  
\[
    A \sum_{n=1} ^N \alpha_n \phi_n =  \sum_{n=1} ^N  -n^2 \alpha_n \phi_n .
\]
  with domain
  \[
     D(A)= \{ x = \sum_{n=1} ^{\infty} \alpha_n \phi_n \in X \mid \sum_{n=1} ^{\infty} |n^2 \alpha_n|^2 < \infty\}
  \]
  We define $C$ as the square root of $-A$, i.e.
  \begin{equation}
    \label{eq:16}
     C \sum_{n=1} ^N \alpha_n \phi_n =  \sum_{n=1} ^N  n \alpha_n \phi_n 
  \end{equation}
  with domain
  \[
     D(C)= \{ x = \sum_{n=1} ^{\infty} \alpha_n \phi_n \in X \mid \sum_{n=1} ^{\infty} |n \alpha_n|^2 < \infty\}.
  \]
  A straightforward calculation gives that for $x_0 = \sum_{n=1} ^N \alpha_n \phi_n$, we have that
  \[
     \int_0^{\infty} \| C T(t) x_0\|^2 dt = \int_0^{N} \sum_{n=1}^{N} | n e^{-n^2t}  \alpha_n|^2 dt = \frac{1}{2} \sum_{n=1}^{N} | \alpha_n|^2 = \frac{1}{2} \|x_0\|^2.
  \]
  Since the finite sums lie dense, we conclude that $C$ is admissible. It is easy to see that $C$ commutes with the semigroup, and thus from Theorem \ref{T2.2} we have that
  \begin{equation}
    \label{eq:17}
    \|C T(t) \|\leq \frac{\gamma}{\sqrt{t}}.
  \end{equation}
  for some $\gamma$ independent of $t$. 

  Next choose $x_0 =\phi_n$ and $t= n^{-2}$. Using (\ref{eq:9}) and (\ref{eq:16}) we see that
  \[
     C T(t) x_0 = n e^{-1} \phi_n = \frac{e^{-1}}{\sqrt{t}} x_0,
  \]
  and thus the estimate (\ref{eq:17}) cannot be improved.
\end{example}

The Lebesgue extension of an admissible operator is defined by
\[
  C_{L} x = \lim_{t\rightarrow 0} \frac{1}{t} C \int_0^{t} T(\tau) x d\tau,
\]
where
\[
  D(C_L) = \{ x \in X \mid \mbox{limit exists} \}.
\]
A similar extension can be define using the resolvent. The Lambda extension of an admissible operator is defined by
\[
  C_{\Lambda} x = \lim_{\lambda \rightarrow \infty} \lambda C (\lambda I - A)^{-1} x,
\]
where
\[
  D(C_\Lambda) = \{ x \in X \mid \mbox{limit exists} \}.
\] 
The relation between these extension is still not completely understood, but for admissible operators which commute with the semigroup, we have that both extensions are closed operators.
\begin{lemma}
\label{L2.7}
  Let $C$ be an admissible operator which commutes with the semigroup, then the same holds for its Lebesgue and Lambda extension. Furthermore, these extensions are closed operators.
\end{lemma}
{\bf Proof}\/ 
\mbox{}
Since $A^{-1}$ and $CA^{-1}$ are bounded, we find for $x_0 \in D(C_L)$
\begin{align*}
   A^{-1} C_L x_0 =&\, A^{-1} \lim_{t \downarrow 0} \frac{1}{t}  C \int_0^t T(\tau) x_0 d \tau = \lim_{t \downarrow 0} \frac{1}{t}   A^{-1} C \int_0^t T(\tau) x_0 d \tau \\
  = &\, \lim_{t \downarrow 0} \frac{1}{t}  C  A^{-1} \int_0^t T(\tau) x_0 d \tau = C  A^{-1} \lim_{t \downarrow 0} \frac{1}{t}  \int_0^t T(\tau) x_0 d \tau \\
  =&\, C A^{-1} x_0= C_L A^{-1} x_0,
\end{align*}
where we have used that $\int_0^t T(\tau) x_0 d \tau \in D(A)$ and $C$ commutes with $A^{-1}$.
This proves the first assertion.

Using once more that $CA^{-1}$ and $A^{-1}$ are bounded, we have for $x_0 \in D(C_L)$
\begin{align*}
  \nonumber
 C A^{-1} \int_0^t T(\tau) x_0 d \tau 
  =&\,  \int_0^t C A^{-1} T(\tau) x_0 d \tau  \\
  \nonumber
  =&\, \int_0^t T(\tau) C A^{-1} x_0 d \tau \\
  =&\, \int_0^t T(\tau) A^{-1} C_L x_0 d \tau  = A^{-1}\int_0^t T(\tau) C_L x_0 d \tau.
\end{align*}
Let $x_n$ be a sequence in $D(C_L)$ which converges to $x \in X$, such that $C_L x_n$ converges to $z \in X$. Then by the above we find that
\begin{equation}
  \label{eq:19}
  C A^{-1} \int_0^t T(\tau) x d \tau = A^{-1}\int_0^t T(\tau) z d \tau
\end{equation}
Since $ \int_0^t T(\tau) x d \tau \in D(A)$, we find that 
\begin{equation}
  \label{eq:20}
  A^{-1}\int_0^t T(\tau) z d \tau =  C A^{-1} \int_0^t T(\tau) x d \tau = A^{-1} C\int_0^t T(\tau) x d \tau.
\end{equation}
Hence we have that
\[
  \int_0^t T(\tau) z d\tau = C \int_0^t T(\tau)x d \tau.
\]
Since $t^{-1} \int_0^t T(\tau) z d\tau$ converges to $z$ for $t\downarrow0$, we conclude from the above equality that $x \in D(C_L)$ and $C_Lx = z$.

The proof for $C_{\Lambda}$ goes very similarly. Basically in the above proof, $\int_{0}^t T(\tau) x d\tau$ is replaced by $(\lambda I - A)^{-1}x$.
\hfill$\Box$
\medskip

By Weiss \cite{Weis94a} we have that $C_{\Lambda}$ is an extension of $C_L$. We claim that for admissible $C$'s which commute with the semigroup they are equal.

\section{${\mathcal H}_{\infty}$-calculus}

For $g \in {\mathcal H}_{\infty}^-$ we define the following mapping from $X$ to $L^2((0,\infty);X)$
\begin{equation}
  \label{eq:10}
  {\mathfrak O}_g x_0 = M_g \left( T(t) x_0 \right).
\end{equation}
Hence we have taken in Definition \ref{D1.1}  $f(t)=T(t)x_0$.

It is clear that ${\mathfrak O}_g$ is a linear bounded operator from $X$ into
$L^2((0,\infty);X)$. Furthermore, from (\ref{eq:3}) we have that
\begin{equation}
  \label{eq:11}
  \sigma_{\tau}\left({\mathfrak O}_g x_0\right) = M_g
  \left(\sigma_\tau\left(T(t)x_0 \right)\right) = M_g T(t+\tau) x_0 =
  {\mathfrak O}_g \left(T(\tau)x_0\right),
\end{equation}
where we have used the semigroup property. Hence ${\mathfrak O}_g$ is
a well-posed output map, and so by Lemma \ref{L:2.1a} we conclude that
${\mathfrak O}_g$ can be written as
\begin{equation}
  \label{eq:12}
  {\mathfrak O}_g x_0 = g(A) T(t) x_0
\end{equation}
for some infinite-time admissible operator $g(A)$ which is bounded from the domain of $A$ to $X$. 

Since for all $t,\tau \in [0,\infty)$ there holds $T(\tau) T(t) = T(t) T(\tau)$, we conclude from (\ref{eq:10}) and (\ref{eq:7}) that
\[
  {\mathfrak O}_g T(t) = T(t) {\mathfrak O}_g, \qquad t \geq 0.
\]
Hence by (\ref{eq:12}), we see that $g(A)$ is an admissible operator
which commutes with the semigroup. Theorem \ref{T2.2} implies that for $t >0$, $g(A)
T(t)$ can be extended to a bounded operator and
\begin{equation}
  \label{eq:3.1}
   \|g(A) T(t) \| \leq \frac{\gamma}{\sqrt{t}}.
\end{equation}
Note that for $t \in [0,1]$ this $\gamma$ can be chosen as $\sup_{t \in [0,1]} \|T(t)\| \cdot \|g\|_{\infty}$. 

The Laplace transform of $ {\mathfrak O}_g$ equals $g(A) (sI-A)^{-1}$. Combining this with the definition of ${\mathfrak O}_g$, implies that
\begin{equation}
  \label{eq:21}
  \|g(A) (sI-A)^{-1}\| \leq \frac{\|g\|_{\infty}}{\sqrt{{\mathrm{Re}}(s)}} \|x_0\|,
\end{equation}
where we have taken the norm in $X$, see also Weiss \cite{Weis91}.

Since we have written this admissible operator as the function $g$
working on the operator $A$, there is likely to be a relation with
functional calculus. This is shown next.
\begin{lemma}
  \label{L3.1} If $g\in {\mathcal H}_{\infty}^-$ is the inverse
  Fourier transform of the function $h$, with $h \in
  L^1(-\infty,\infty)$ with support in $(-\infty,0)$, then $g(A)$ is
  bounded
  \begin{equation}
  \label{eq:1.5}
    g(A)x_0= \int_0^{\infty} T(t) h(-t)x_0 dt,
  \end{equation}
  and so $g(A)$ corresponds to the classical definition of the
  function of an operator.
\end{lemma}

So if $g$ is the Fourier transform of an absolutely integrable function, then $g(A)$ is bounded. We would like to know when it is bounded for every $g$. For this, we extend the definition of ${\mathfrak O}_g$.

Let $C$ be an admissible output operator for the semigroup
$\left(T(t)\right)_{t\geq 0}$. By definition, we know that $CT(\cdot) x_0
\in L^2((0,\infty);Y)$ for all $x_0 \in X$. We define
\begin{equation}
  \label{eq:13}
  \left(C \circ {\mathfrak O}_g\right) x_0 = M_g \left( CT(t) x_0 \right)
\end{equation}
It is clear that this is a bounded mapping from $X$ to
$L^2((0,\infty);Y)$.

As before we have that
\begin{equation}
  \label{eq:14}
  \sigma_{\tau} \left(\left(C \circ {\mathfrak O}_g\right)( x_0)\right) = \left(C \circ {\mathfrak O}_g\right)\left(T(\tau)x_0\right).
\end{equation}
And so we can write $\left(C \circ {\mathfrak O}_g\right) x_0$ as
$\tilde{C}_g T(\cdot) x_0$ for some infinite-time admissible
$\tilde{C}_g$. We have that
\begin{lemma}
  \label{L3.2}
  The infinite-time admissible operator $\tilde{C}_g$ satisfies
\begin{equation}
\label{eq:15}
   \tilde{C}_gx_0 = C g(A) x_0, \qquad \mbox{for } x_0 \in D(A^2).
\end{equation}
\end{lemma}
{\bf Proof} For $x_0 \in D(A^2)$, we introduce $x_1 = A x_0$. Then the following equalities hold in $L^2((0,\infty);Y)$.
\begin{eqnarray*}
  \tilde{C}_g T(t) x_0 &=& \left(C \circ {\mathfrak O}_g\right) x_0\\
  &=&
  M_g \left( CT(t) x_0 \right) \\
  &=&
   M_g \left( CT(t) A^{-1} x_1 \right) \\
  &=&
  M_g \left( CA^{-1} T(t) x_1 \right)\\
  &=&
  CA^{-1} g(A) T(t) x_1\\
  &=& C g(A) T(t) A^{-1} x_1 = C g(A) T(t) x_0,
  \end{eqnarray*}
  where we have used (\ref{eq:7}).
  Since both functions are continuous at zero, we find that (\ref{eq:15}) holds.
\hfill$\Box$
\medskip

Based on this result, we denote $\tilde{C}_g$ by $C \circ g(A)$. 

Using this, we can prove the following theorems.
\begin{theorem}
  \label{T3.2a}%
  The mapping $g \mapsto g(A)$ forms a (unbounded) ${\mathcal
    H}_{\infty}^-$-calculus.
\end{theorem}
{\bf Proof}\/ It only remains to show that $(g_1g_2)(A)= g_1(A) g_2(A)$. By Lemma \ref{L1.2} we have that
\[
 {\mathfrak O}_{g_1g_2}x_0 = M_{g_1g_2}\left(T(t) x_0\right) = M_{g_1}M_{g_2}\left(T(t) x_0\right).
\]
For $x_0 \in D(A)$ the last expression equals $M_{g_1} \left( g_2(A) T(t) x_0 \right)$, see (\ref{eq:12}). Since $g_2(A)$ commutes with the semigroup, we find that
\[
  {\mathfrak O}_{g_1g_2}x_0 = M_{g_1} \left(  T(t) g_2(A) x_0 \right).
\]
Using (\ref{eq:12}) twice, we obtain 
\[
   (g_1g_2)(A) T(t) x_0 = {\mathfrak O}_{g_1g_2}x_0= g_1(A) T(t) g_2(A) x_0
\]
This is an equality in $L^2((0,\infty);X)$. However, if we take $x_0 \in D(A^2)$, then this holds point-wise, and so for $x_0 \in D(A^2)$.
\[
  (g_1g_2)(A)x_0 = g_1(A) g_2(A) x_0
\]
This concludes the proof.
\hfill$\Box$

\begin{theorem}
  \label{T3.3}
  If there exists an admissible $C$ such that $(C,A)$ is exactly observable, i.e., these exists an $m_1>0$ such that for all $x_0 \in X$ there holds
  \[
    \int_0^{\infty} \|C T(t) x_0 \|^2 dt \geq m_1 \|x_0\|^2
  \]
  then $g(A)$ is bounded for every $g \in {\mathcal H}_{\infty}^-$.
  Furthermore, if $m_2$ is the admissibility constant, see equation (\ref{eq:2.1}), then
  \begin{equation}
    \label{eq:4}
    \|g(A) \| \leq \sqrt{\frac{m_2}{m_1}} \|g\|_{\infty}.
  \end{equation}
\end{theorem}
{\bf Proof}\/
Let $x_0 \in D(A^2)$
\begin{eqnarray*}
  m_1 \|g(A)x_0\|^2
  &\leq&  \|C T(t) g(A)x_0\|^2_{L^2((0,\infty);Y)}\\
  &=&
   \|Cg(A) T(t) x_0\|^2_{L^2((0,\infty);Y)}\\
  &=&
   \|C \circ {\mathfrak O}_g x_0 \|^2_{L^2((0,\infty);Y)}\\
  &\leq& \|g\|_{\infty}^2 \|C T(t) x_0\|^2_{L^2((0,\infty);Y)}\\
 &\leq&
   m_2 \|g\|_{\infty}^2 \|x_0\|^2.
\end{eqnarray*}
Since $D(A^2)$ is dense, we obtain the result.
\hfill$\Box$
\medskip

As a corollary we obtain the well-known von Neumann inequality. Recall that the operator $A$ is dissipative if
\begin{equation}
  \label{eq:25}
  \langle x_0, Ax_0 \rangle + \langle  Ax_0, x_0\rangle \leq 0 \qquad \mbox{ for all } x_0 \in D(A).
\end{equation}
\begin{corollary}
  \label{C3.3.a} If $A$ is a dissipative operator and its corresponding semigroup is exponentially stable, then $A$ has a bounded ${\mathcal H}_{\infty}^-$ calculus and for all $g \in {\mathcal H}_{\infty}^-$
  \begin{equation}
    \label{eq:5}
    \|g(A) \| \leq \|g\|_{\infty}.
  \end{equation}
\end{corollary}
{\bf Proof}\/
  Since $A$ is dissipative and since its semigroups is exponentially stable, we have that $A^{-1}$ is bounded and dissipative. We define $Q$ via
  \begin{equation}
    \label{eq:6}
    \langle x_1,Qx_2\rangle =-\langle  A^{-1}x_1, x_2\rangle - \langle x_1, A^{-1}x_2 \rangle , \qquad x_1,x_2 \in X.
  \end{equation}
  It is easy to see that $Q$ is bounded, self-adjoint and by the dissipativity of $A^{-1}$ we have that $Q \geq 0$. Define on the domain of $A$ the operator $C$ as $C=\sqrt{Q}A$, then from (\ref{eq:6}) we find that
  \begin{equation}
    \label{eq:18}
    - \langle C x_1,C x_2\rangle = \langle x_1, Ax_2 \rangle + \langle  Ax_1, x_2\rangle, \qquad x_1,x_2 \in D(A).
  \end{equation}
   Combining this Lyapunov equation with the exponential stability, gives that for all $x_0 \in D(A)$
  \begin{equation}
    \label{eq:8}
    \int_0^{\infty} \|C T(t)x_0 \|^2dt = \|x_0\|^2.
  \end{equation}
  Thus we see that the constants $m_1$ and $m_2$ in Theorem \ref{T3.3} can be chosen to be one, and so (\ref{eq:4}) gives the results.
\hfill$\Box$
\medskip

If $A$ generates an exponentially stable semigroup and if there exists an admissible $C$ for which $(C,A)$ is exactly observable, then it is not hard to show that the semigroup is similar to a contraction semigroup. Using this, one can also obtain the above result by Theorem G of \cite{AlDM96}. The following result has been proved by McIntosh in \cite{McIn86}.
\begin{theorem}
  \label{T3.4} Assume that $A$ generates an exponentially stable semigroup. If $(-A)^{\frac{1}{2}}$ is admissible for $\left(T(t)\right)_{t \geq 0}$ and  $\left(-A^*\right)^{\frac{1}{2}}$ is admissible for the adjoint semigroup $\left( T(t)^* \right)_{t \geq 0}$, then $g(A)$ is bounded for every $g \in {\mathcal H}_{\infty}^-$. Thus this semigroup has a bounded ${\mathcal H}_{\infty}^-$-calculus.
\end{theorem}
{\bf Proof}\/
Since $A^{1/2}$ is admissible, Lemma \ref{L3.2} gives that $A^{1/2}
\circ g(A)$ is also admissible.
Consider for $x_1 \in D(A^*)$ and $x_0 \in D(A^2)$ the following
\begin{align*}
   \langle x_1, g(A)& x_0 \rangle - \langle x_1, g(A) T(t) x_0 \rangle \\
  & = 
  \int_0^t \langle x_1, (-A) T(\tau)  g(A) x_0 \rangle d \tau \\
  &=
  \int_0^t \langle \left(-A^*\right)^{\frac{1}{2}} x_1,  \left(-A\right)^{\frac{1}{2}} g(A)  T(\tau) x_0 \rangle d \tau \\
  &=
   \int_0^t \langle \left(-A^*\right)^{\frac{1}{2}} T(\frac{\tau}{2})^* x_1,   g(A) \left(-A\right)^{\frac{1}{2}}  T(\frac{\tau}{2}) x_0 \rangle d \tau \\
  &\leq \sqrt{\int_0^t \| \left(-A^*\right)^{\frac{1}{2}} T(\frac{\tau}{2})^* x_1 \|^2 d\tau} \sqrt{\int_0^t \|g(A) \left(-A\right)^{\frac{1}{2}}  T(\frac{\tau}{2}) x_0 \|^2 d\tau}\\
  &\leq \sqrt{\int_0^t \| \left(-A^*\right)^{\frac{1}{2}} T(\frac{\tau}{2})^* x_1 \|^2 d\tau} \|g\|_{\infty} \sqrt{\int_0^{\infty} \| \left(-A\right)^{\frac{1}{2}} T(\frac{\tau}{2}) x_0 \|^2 d\tau}\\
  &\leq  m_1 \|x_1\| m_2\|g\|_{\infty} \|x_0\|,
\end{align*}
where $m_1$ and $m_0$ are the admissibility constant of $ \left(-A^*\right)^{\frac{1}{2}} $ and $ \left(-A^*\right)^{\frac{1}{2}} $, respectively. Furthermore, we used (\ref{eq:1.1}).

Since the sets $D(A^*)$ and $D(A^2)$ are dense in $X$, we obtain that
\begin{equation}
  \label{eq:22}
 \|g(A) \|  \leq m_1 m_2 \|g\|_{\infty} + \|g(A) T(t)\|.
\end{equation}
By Theorem \ref{T2.2} we know that $g(A) T(t)$ is bounded, and so we conclude that $\left(T(t)\right)_{t \geq 0}$ has a bounded  ${\mathcal H}_{\infty}^-$-calculus. 
\hfill$\Box$
\medskip

In McIntosh \cite{McIn86} the above  theorem was proved using square function estimates. The admissibility of $(-A)^{\frac{1}{2}}$ can be written as
\begin{eqnarray*}
   m \|x_0 \|^2 &\geq& \int_0^{\infty} \|(-A)^{\frac{1}{2}} T(t) x_0 \|^2 dt \\
   &=& \int_0^{\infty}  \|(-tA)^{\frac{1}{2}} T(t) x_0 \|^2 \frac{dt}{t}.
\end{eqnarray*}
The latter is the ``square function estimate'' for $\psi(s) = (-s)^{\frac{1}{2}} e^{s}$, and so the admissibility condition can be seen as a square function estimate. The other condition used in \cite{McIn86} is that the operator $A$ is sectorial on a sector larger than the sector on which the scalar functions are defined. Since we have as function class ${\mathcal H}_{\infty}^-$ and since our operators $A$ are assumed to generate an exponential semigroup, this condition seems not to satisfied. However, the admissibility assumptions made in the theorem imply that $A$ generates a bounded analytic semigroup, and so the condition of McIntosh is satisfied.
\begin{lemma}
  Let $A$ generate an exponentially stable semigroup and let $(-A)^{\frac{1}{2}}$  and  $\left(-A^*\right)^{\frac{1}{2}}$ be admissible operators for for $\left(T(t)\right)_{t \geq 0}$ and $\left( T(t)^* \right)_{t \geq 0}$, respectively. Then $A$ generates a bounded analytic semigroup.
\end{lemma}
{\bf Proof}\/ The proof is similar to the proof of Theorem \ref{T2.2}. Let $x_1 \in D(A^*)$ and $x_0 \in D(A)$. Then for $t>0$ we find
\begin{align*}
    t \langle x_1, AT(t) x_0 \rangle 
  & = 
  \int_0^t \langle x_1, A T(t) x_0 \rangle d \tau \\
  &=
  -\int_0^t \langle \left(-A^*\right)^{\frac{1}{2}} x_1,  \left(-A\right)^{\frac{1}{2}}  T(t) x_0 \rangle d \tau \\
  &=
   -\int_0^t \langle \left(-A^*\right)^{\frac{1}{2}} T(\tau)^* x_1, \left(-A\right)^{\frac{1}{2}}  T(t-\tau) x_0 \rangle d \tau \\
  &\leq \sqrt{\int_0^t \| \left(-A^*\right)^{\frac{1}{2}} T(\tau)^* x_1 \|^2 d\tau} \sqrt{\int_0^t \| \left(-A\right)^{\frac{1}{2}}  T(t-\tau ) x_0 \|^2 d\tau}\\
  &\leq  m_1 \|x_1\| m_2\|x_0\|,
\end{align*}
where we used that $(-A)^{\frac{1}{2}}$  and  $\left(-A^*\right)^{\frac{1}{2}}$ are admissible. Since the domain of $A^*$ and $A$ are dense, we obtain that
\[
   \|AT(t) \| \leq \frac{M}{t}, \qquad t >0
\]
By Theorem II.4.6 of \cite{EnNa00}, we conclude that generates a bounded analytic semigroup.
\hfill$\Box$
\medskip

\mbox{}From \cite{McIn86} we know that if the conditions of Theorem \ref{T3.4} hold, then is the semigroup similar to a contraction (or $(-A)^{\frac{1}{2}}$ is exactly observable). We show this next.
\begin{lemma}
Under the condition of Theorem \ref{T3.4} we have that $(-A)^{\frac{1}{2}}$ is exactly observable, and thus $\left( T(t)\right)_{t\geq 0}$ is similar to a contraction.
\end{lemma}
{\bf Proof}\/
In idea the proof is the same as that of Theorem \ref{T3.4}. Let $x_1\in D(A^*)$ and $x_0 \in D(A)$ We have that
 \begin{align}
  \nonumber
 \langle x_1, x_0 \rangle 
  & = 
  \int_0^{\infty} \langle x_1, (-A) T(\tau) x_0 \rangle d \tau \\
  \nonumber
  &=
  \int_0^{\infty} \langle \left(-A^*\right)^{\frac{1}{2}} x_1,  \left(-A\right)^{\frac{1}{2}}  T(\tau) x_0 \rangle d \tau \\
  \label{eq:24}
    &=
   \int_0^{\infty} \langle \left(-A^*\right)^{\frac{1}{2}} T(\frac{\tau}{2})^* x_1,  \left(-A\right)^{\frac{1}{2}}  T(\frac{\tau}{2}) x_0 \rangle d \tau.
\end{align}
Hence 
\begin{align*}
   |\langle x_1, x_0 \rangle |  &\leq \sqrt{\int_0^\infty \| \left(-A^*\right)^{\frac{1}{2}} T(\frac{\tau}{2})^* x_1 \|^2 d\tau} \sqrt{\int_0^\infty \| \left(-A\right)^{\frac{1}{2}}  T(\frac{\tau}{2} ) x_0 \|^2 d\tau}\\
  &\leq  m_1 \|x_1\|\sqrt{\int_0^\infty \| \left(-A\right)^{\frac{1}{2}}  T(\frac{\tau}{2} ) x_0 \|^2 d\tau}
\end{align*}
Since the domain of $A^*$ is dense we conclude that
\begin{equation}
  \label{eq:26}
  \|x_0\| = \sup_{x_1 \neq 0}  \frac{|\langle x_1, x_0 \rangle | }{\|x_1\|} \leq m_1 \sqrt{\int_0^\infty \| \left(-A\right)^{\frac{1}{2}}  T(\frac{\tau}{2} ) x_0 \|^2 d\tau}.
\end{equation}
Thus $(-A)^{\frac{1}{2}}$ is exactly observable.
\hfill$\Box$
\medskip

We remark that with the above result, Theorem \ref{T3.4} follows also from Theorem \ref{T3.3}. However, we decided to present this independent proof.

\subsection*{Acknowledgement}

The author want to thank Markus Haase, Bernhard Haak, and Christian Le Merdy whom have helped him to understand functional calculus.

\end{document}